\documentclass{article}
\usepackage{color}
\usepackage{graphicx}
\usepackage{amsmath}
\usepackage{amssymb}
\usepackage{mathrsfs}
\usepackage[numbers,sort&compress]{natbib}
\usepackage[colorlinks,linkcolor=black,anchorcolor=black,citecolor=black]{hyperref}
\usepackage{amsthm}
\usepackage{pgf}
\usepackage{graphicx}
\usepackage{tikz}

\usepackage{graphicx}
\usepackage[T1]{fontenc}
\usepackage[english]{babel}
\usepackage{listings}
\usepackage{xcolor,mathrsfs,url}
\usepackage{amssymb}
\usepackage{amsmath}
\usepackage{ifthen}
\usepackage{caption}
\usepackage{float}

\usepackage{authblk}

{\LARGE }

\newcommand{\AAA}{{\mathcal A}}
\newcommand{\ddddd}{\mathrm{d}}

\newcommand{\deff}{\overset{\mathrm{def}}{=}}
\newcommand{\yl}{\left[\begin{matrix}
		1\\0
	\end{matrix}\right]}

\begin{document}
	\title{$L^2$-Sobolev space bijectivity of the scattering-inverse scattering transforms related to defocusing Ablowitz-Ladik systems}
	\author[a,b]{Meisen Chen}
	\author[c]{Engui Fan}
	\author[a]{Jingsong He\thanks{Corresponding author: hejingsong@szu.edu.cn}}
	\affil[a]{\small Institute for Advanced Study, Shenzhen University, Shenzhen, 518060, China}
	\affil[b]{College of Physics and Optoelectronic Engineering, Shenzhen University, Shenzhen, 518060, China}
	\affil[c]{School of Mathematical Sciences, Fudan University, Shanghai, 200433, China}
	\maketitle
	
	\theoremstyle{plain}
	\newtheorem{proposition}{Proposition}[section]
	\newtheorem{lemma}[proposition]{Lemma}
	\newtheorem{theorem}[proposition]{Theorem}
	\newtheorem{drhp}[proposition]{$\bar\partial$-RH problem}
	\newtheorem{rhp}[proposition]{RH problem}
	\newtheorem{dbarproblem}[proposition]{$\bar\partial$-problem}
	\newtheorem{remark}[proposition]{Remark}

\baselineskip=20 pt

\begin{abstract}
	In this paper, we establish $L^2$-Sobolev space bijectivity of the inverse scattering transform related to the defocusing Ablowitz-Ladik system. 
	On the one hand, in the direct problem, based on the spectral problem, we establish the reflection coefficient and  the corespondent Riemann-Hilbert problem. 
	And we also prove that if the potential belongs to $l^{2,k}$ space, then the reflection coefficient belongs to $H^k_\theta(\Sigma)$.
	On the other hand, in the inverse problem, based on the Riemann-Hilbert problem, we obtain the corespondent reconstructed formula and recover potentials from reflection coefficients. 
	And we also confirm that if reflection coefficients are in $H^k_\theta(\Sigma)$, then we show that potentials also belong to $l^{2,k}$. 
	This study also confirm that for the initial-valued problem of defocusing Ablowitz-Ladik equations, it the initial potential belongs to $l^{2,k}$ and satisfying $\parallel q\parallel_\infty<1$, then the solution for $t\ne0$ also belongs to $l^{2,k}$.
	\\[5pt]
\end{abstract}

\newpage
\tableofcontents

\section{Introduction}
\indent

In this paper, we consider the discrete spectral problem
\begin{align}\label{e2}
	V(z,n+1)=(z^{\sigma_3}+Q_n)V(z,n),
\end{align}
where $n$ is the discrete spatial variable, $z$ is the spectral variable, and 
\begin{align*}
	\sigma_3=\left[\begin{matrix}
		1&0\\0&-1
	\end{matrix}\right],\quad
	Q_n=\left[\begin{matrix}
		0&q_n\\\bar q_n&0
	\end{matrix}\right]. 
\end{align*}
The spectral problem (\ref{e2}) is firstly introduced by Ablowitz and Ladik \cite{Ablowitz1975nonlinear} in 1975, and associated to the defocusing Ablowitz-Ladik system
\begin{align}\label{e1}
	i\partial_tq_n=q_{n+1}-2q_n+q_{n-1}-|q_n|^2(q_{n+1}+q_{n-1}),
\end{align}
which is an integrable difference-differential equation and the discretization of nonlinear Schr\"odinger equation
\begin{align*}
	i\partial_tu=\partial^2_xu-|u|^2u.
\end{align*} 
In what follows, we properly denote the discrete potential $\{q_n\}_{n=-\infty}^\infty$ as $q$ without confusion of notation. 
In this paper, the purpose is to investigate inverse scattering mapping for (\ref{e2}) and potentials satisfying that
\begin{align}\label{e3s}
	\parallel q\parallel_\infty< 1\quad\text{and}\quad  q\in l^{2,k},
\end{align} 
where $l^{2,k}$ denote the discrete weighted Sobolev space
\begin{align}
	l^{2,k}=\left\{\{q_n\}_{n=-\infty}^\infty:\sum_{n=-\infty}^{\infty}(1+n^2)^k|q_n|^2<\infty\right\},\quad k>0,
\end{align}
and we denote $\parallel q\parallel_{2,k}$ as the $l^{2,k}$-norm of $q$:
\begin{align*}
	\parallel q\parallel_{2,k}=\left(\sum_{n=-\infty}^{\infty}(1+n^2)^k|q_k|^2\right)^\frac{1}{2}.
\end{align*} 

In this paragraph, we start the main result in this paper. 
In the direct problem, with potentials in the discrete weighted Sobolev space, we rigorously prove that the reflection coefficient belongs to
\begin{align*}
	H^k_\theta=H^k_\theta(\Sigma)=\left\{f(\theta)\in L^2_\theta(\Sigma):\partial_\theta^\alpha f\in L^2_\theta(\Sigma),\alpha=1,\dots,k\right\},
\end{align*}
where
\begin{align*}
	\Sigma=\left\{z=e^{i\theta}:\theta\in[-\pi,\pi]\right\}
\end{align*}
is the jump contour.  
By Fourier analysis, we also denote the $H^k_\theta$-norm of $f\equiv f(\theta)$ as
\begin{align*}
	\parallel f\parallel_{H^k_\theta}=\left(\sum_{n=-\infty}^{\infty}(1+n^2)^k|\hat f(n)|^2\right)^\frac{1}{2}, \quad \hat f(n)=\frac{1}{2\pi}\int_{-\pi}^{\pi}f(\theta)e^{-in\theta}\ddddd\theta.
\end{align*}
where $\hat f(n)$ denote the $n$th entry in Fourier series of $f(\theta)$. 
In the inverse problem, we also prove that when the reflection coefficient belongs to $H^k_\theta$, the potential is also proven in $l^{2,k}$. 
We see that this problem is basically settled from the point of view of Fourier theory. 
Setting 
\begin{align}\label{e5s}
	c_{-\infty}=\prod_{n=-\infty}^{\infty}(1-|q_n|^2), 
\end{align}
we can see  from (\ref{e1}) that $c_{-\infty}$ is time-independent by basic computation.
In the inverse problem, we also notice that the constant $c_{-\infty}$ can be equivalently represented by the reflection coefficient:
\begin{align}\label{e6}
	c_{-\infty}=\exp\left[\frac{1}{2\pi}\int_{-\pi}^{\pi}\ln(1-|r(\theta)|^2)\ddddd\theta\right].
\end{align}

\begin{figure}
	\centering
	\includegraphics[width=0.8\linewidth]{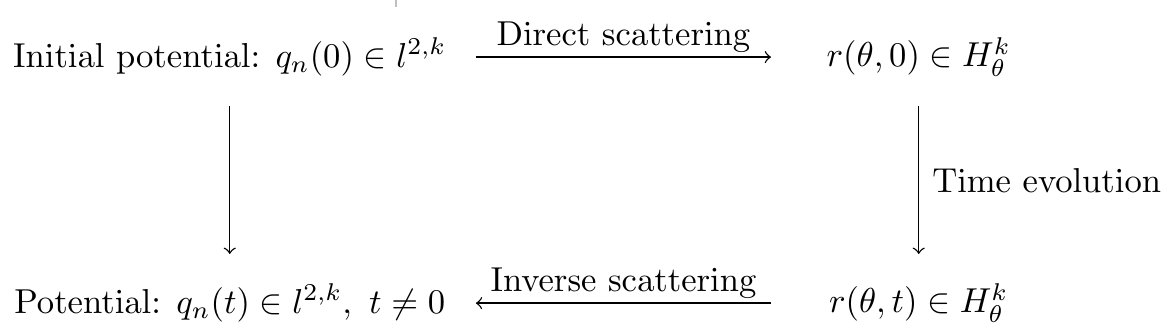}
	\caption{Time evolution for defocusing Ablowitz-Ladik lattice}\label{f1}
\end{figure}

At the end of this paper, we discuss the time evolution of solutions for defocusing Ablowitz-Ladik systems. 
When the initial potential $q(0)$ belongs to $l^{2,k}$ and satisfies that $\parallel q(0)\parallel_\infty<1$, the reflection potential $r(\theta,0)$ belongs to $H^k_\theta$ as shown in Section \ref{s2}, and the flow $t\rightarrow r(\theta,t)$ persists reflection coefficients in $H^k_\theta$.
Then, as shown in Section \ref{s3}, for $t\ne0$, the initial-valued problem for (\ref{e1}) has solution and it also belongs to $l^{2,k}$. See Figure \ref{f1}. 

The inverse scattering transform is an important method for integrable systems and has made a great progress.
In 1967, the inverse scattering transform was firstly introduced when solving the KdV equation by Gardne, Greene, Kruskal and Miura \cite{gardner1967method}.
In 1974, Ablowitz et al solved ZS-AKNS systems by inverse scattering transforms \cite{Ablowitz1974the}. 
In 1975, Shabat \cite{shabat1975inverse} investigated the inverse scattering tranvsform by Riemann-Hilbert (RH) method. 
Since 1980s, the RH method has been applied to many integrable systems \cite{Wang2010integrable,Biondini2014inverse,Biondini2015inverse,Kraus2015the,Pichler2017on,xu2014a,xiao2016a,kang2018multi,yang2019high,yang2010nonlinear,yang2021riemann}.
Also, it was generalized to solve the discrete integrable systems \cite{Ablowitz2007inverse,ablowitz2020discrete,Prinari,Ortiz2019inverse,chen2021riemann}. 
Based on the RH method, in 1993, Deift et al \cite{deift1993a,deift1994long} developed a nonlinear steepest descent method, well-known as Deift-Zhou Method, to investigate the long-time behavior of defocusing NLS equations and MKdV equations with initial potentials in Schwartz space. 
Since then, people have applied Deift-Zhou method to long-time asymptotic analysis of many integrable systems \cite{grunert2009long,de2009long,xu2015long,xu2018long,huang2015long,zhu2018the}. 
Researchers also generalized Deift-Zhou method to apply it on discrete integrable systems with Schwartz initial potentials, such as Toda lattice \cite{kruger2009long,kruger2009long2}, discrete nonlinear Schr\"odinger equations \cite{yamane20114long,yamane2015long,yamane2019long} and defocusing discrete mKdV equations \cite{chen2020long}. 

Since these progresses of the inverse scattering transform, researchers concerned how the inverse spectral mapping works from the potential to the reflection coefficient.
For AKNS systems, Beals and Coifman \cite{beals1984scattering} investigated the inverse scattering transform and proved that the potential belongs to Schwartz space if and only if the reflection coefficient is in Schwartz space. 
In 1998, Zhou \cite{zhou1998l2} established the $L^2$-Sobolev space bijectivity for the inverse spectral transform associated with ZS-AKNS systems including nonlinear Schr\"odinger equations. 
In \cite{zhou1998l2}, it was proven that if the potential belongs to $L^2((1+x^2)^k\ddddd x)\cap H^l$, for $k\ge1$ and $l\ge0$, then the reflection coefficient belongs to $H^k\cap L^2((1+x^2)^l\ddddd x)$, and vice versa.
In 2019, Liu \cite{Liu2019l2} generalized the result in \cite{zhou1998l2} and established $L^2$-Sobolev space bijectivity of the inverse scattering transform associated with $3\times3$ AKNS systems. 
These studies are valuable for long-time asymptotic analysis for integrable systems of initial potential with lower smoothness. 
In 2003, based on \cite{zhou1998l2}, Deift and Zhou \cite{deift2003long} analyzed the long-time behavior of NLS equation with initial data in weighted Sobolev space. 
In 2018, Borghese et al \cite{Borghese2018long} investigate long-time asymptotics for focusing NLS equations with soliton and initial data in weighted Sobolev space by dbar steepest descent method. 
Because of these valuable results for the continuous integrable system, we want to study the generalized result about the spectral problem associated to defocusing Ablowitz-Ladik systems, which we have stated in the first paragraph. 
This study is also valuable for further analyzing properties of the solution for discrete integrable systems with spectral problem (\ref{e2}) and potentials in $l^{2,k}$, for instance, the long-time asymptotic analysis. 
Our study is also the generalization of inverse spectral transform associated with (\ref{e1}) shown in \cite{Ablowitz1975nonlinear,Ablowitz1976nonlinear,ablowitz2004discrete}, in which the potential belongs to the discrete Schwartz space
\begin{align*}
	\left\{\{q_n\}_{n=-\infty}^\infty:\sum_{n=-\infty}^{\infty}(1+n^2)^k|q_n|^2<\infty,\ k\ge1\right\}.
\end{align*}

This article is organized as follows. In section \ref{s2}, we investigate the direct scattering part and establish the mapping from potentials to reflection coefficients. 
In section \ref{s3}, we recover the potential by the reconstructed formula and explain the relation between $q$ and $r$.
In the final section, we discuss the time evolution of the reflection coefficient for defocusing Ablowitz-Ladik systems.

\section{Direct scattering problem}\label{s2}
\indent

In this section, we introduce the direct scattering problem, and the purpose is to construct the reflection coefficient and RH problem \ref{r3} associated to it,  and then prove Theorem \ref{t2.4}.

\subsection{Jost solutions, scattering matrix, reflection coefficients and RH problem}
\noindent

Here, we construct the Jost solution, the scattering matrix, the reflection coefficient and the RH problem. 
Since the potential belongs to $l^{2,k}$, $k\ge1$, $q_n$ decays to $0$ as $n\to\pm\infty$ and the spectral problem (\ref{e2}) admits the Jost solution: $X^{(\pm)}=X^{(\pm)}(z,n)$, such that as $n\to\pm\infty$,
\begin{align}\label{e3}
	X^{(\pm)}(z,n)\sim z^{n\sigma_3}.
\end{align}
Naturally introduce the modified Jost solution:
\begin{align}\label{e4}
	Y^{(\pm)}=Y^{(\pm)}(z,n)=z^{-n\sigma_3}X^{(\pm)}(z,n),
\end{align}
which admits the property that as $n\to\pm\infty$,
\begin{align*}
	Y^{(\pm)}(z,n)\sim I. 
\end{align*}
From the spectral problem (\ref{e2}), $\parallel q\parallel_\infty<1$ and (\ref{e3}), we learn that
\begin{align*}
	c_n\deff(\det X^{(+)}(z,n))^{-1}=\prod_{k=n}^\infty(1-|q_n|^2)>0.
\end{align*}
By basic computation, (\ref{e2}), (\ref{e3}) and (\ref{e4}), $Y^{(\pm)}$ admits the following summation equations:
\begin{subequations}\label{e5}
	\begin{align}
		Y^{(+)}(z,n)&=I-z^{-\sigma_3}\sum_{k=n}^{\infty}z^{-k\hat\sigma_3}Q_kY^{(+)}(z,k),\\
		Y^{(-)}(z,n)&=I+z^{-\sigma_3}\sum_{k=-\infty}^{n-1}z^{-k\hat\sigma_3}Q_kY^{(-)}(z,k).
	\end{align}
\end{subequations}
Recalling that  $q\in l^{2,k},\ k\ge1$,
by Swartz inequality, it also satisfies that
\begin{align*} \sum_{n=-\infty}^{\infty}(1+n^2)^\frac{1}{2}|q_n|<\infty.
\end{align*}
Then, by taking Neumann series of $Y^{(\pm)}$ in the summation equation (\ref{e5}), it follows the analyticity for $Y^{(\pm)}=(Y^{(\pm)}_1,Y^{(\pm)}_2)$, which is shown in Lemma \ref{l1}. 
\begin{lemma}\label{l1}
	If $q\in l^{2,k}$, $k\ge1$, then $Y_1^{(-)}$ and $Y_2^{(+)}$ are analytic on $D_+=\{|z|>1\}$, and continuously extended to $D_+\cup\Sigma$.
	In the meanwhile, $Y_2^{(-)}$ and $Y_1^{(+)}$ are analytic on $D_-=\{|z|<1\}$, and continuously extended to $D_-\cup\Sigma$. 
\end{lemma}
\noindent
Because of that $\parallel q\parallel_\infty<1$, by (\ref{e2}) and (\ref{e3}), we always have
\begin{align*}
	\det X^{(+)}\ne 0,\quad \det X^{(-)}\ne0,
\end{align*}
which means that $X^{(+)}$ and $X^{(-)}$ are both invertible eigenfunctions for (\ref{e2}); therefore, there is a matrix $S(z)$, well-known as the scattering matrix and independent on $n$, such that
\begin{align}\label{e8s}
	X^{(-)}(z,n)=X^{(+)}(z,n)S(z),\quad S(z)=\left[\begin{matrix}
		a(z)&\breve{b}(z)\\b(z)&\breve{a}(z)
	\end{matrix}\right]. 
\end{align} 
Since $X^{(\pm)}(z,n)$ is the Jost solution of (\ref{e1}), it's readily seen that $\sigma_1\overline{X^{(\pm)}(\bar z^{-1},n)}\sigma_1$ is a solution of (\ref{e2}) and as $n\to\pm\infty$,
\begin{align*}
	\sigma_1\overline{X^{(\pm)}(\bar z^{-1},n)}\sigma_1\sim z^{n\sigma_3},
\end{align*}
where $\sigma_1=\left[\begin{matrix}
	0&1\\1&0
\end{matrix}\right]$; then, by uniqueness of solution, it follows that
\begin{align*}
	X^{(\pm)}(z,n)=\sigma_1\overline{X^{(\pm)}(\bar z^{-1},n)}\sigma_1,
\end{align*}
which combined with (\ref{e8s}) derive that
\begin{align}\label{e10ss}
	\breve a(z)=\overline{a(\bar z^{-1})},\quad \breve b(z)=\overline{b(\bar z^{-1})}.
\end{align}
By the spectral problem (\ref{e2}), (\ref{e8s}) and (\ref{e10ss}), it's readily seen that
\begin{align}\label{e11s}
\det S(z)=|a(z)|^2-|b(z)|^2=c_{-\infty}>0,\quad z\in\Sigma.
\end{align}
By (\ref{e4}), (\ref{e8s}) and Cramer's rule, we have
	\begin{align}\label{e8}
		&a(z)=c_n\det[Y_1^{(-)}(z,n),Y_2^{(+)}(z,n)],\quad b(z)=c_n\det[Y_1^{(+)}(z,n),Y_1^{(-)}(z,n)],
	\end{align}
which, by taking $n\to+\infty$, is equivalent to
	\begin{align*}
		a(z)&=1+\sum_{k=-\infty}^{\infty}z^{-(2k+1)}q_kY^{(-)}_{21}(z,k),\quad
		b(z)=\sum_{k=-\infty}^{\infty}z^{2k+1}\bar q_kY^{(-)}_{11}(z,k).
	\end{align*}
By (\ref{e8}) and Lemma \ref{l1}, it follows that $a(z)$ is analytic on $D_-$ and continuously extended to $D_-\cup\Sigma$; in the meanwhile, $b(z)$ is continuous on $\Sigma$. 
By WKB expansion, we can derive that $[Y^{(+)}_1,Y^{(-)}_2]$ and $[Y^{(-)}_1,Y^{(+)}_2]$ admit the following asymptotic properties at $z=0$ and $z\to\infty$, respectively:
\begin{subequations}\label{e10s}
	\begin{align}
		&[Y_1^{(+)},Y_2^{(-)}](z,n)\sim z^{-n\hat\sigma_3}\left(\left[\begin{matrix}
			c_n^{-1}&q_{n-1}z\\-c_n^{-1}\bar q_nz&1
		\end{matrix}\right]+\mathcal{O}(z^2)\right),\quad z\to0,\label{e10sa}\\
		&[Y_1^{(-)},Y_2^{(+)}](z,n)\sim z^{-n\hat{\sigma}_3}\left(\left[\begin{matrix}
			1&-c_n^{-1}\frac{q_n}{z}\\\frac{\bar q_{n-1}}{z}&c_n^{-1}
		\end{matrix}\right]+\mathcal{O}\left(\frac{1}{z^2}\right)\right),\quad z\to\infty.\label{e10sb}
	\end{align}
\end{subequations}
From (\ref{e8}) and (\ref{e10sb}), we see that
	\begin{align}\label{e9ss}
		&a(z)\sim 1+\mathcal{O}(z^{-1}),\quad z\to\infty.
	\end{align}
Introduce reflection coefficients
\begin{align*}
	r(z)=b(z)/a(z), \quad z\in\Sigma,
\end{align*}
and a holomorphic function on $\mathbb{C}\setminus \Sigma$:
\begin{align}\label{e9s}
	M\equiv M(z,n)=\begin{cases}
		z^{n\hat\sigma_3}\left(\left[\begin{matrix}
			1&0\\0&c_n
		\end{matrix}\right]\left[\frac{Y_1^{(-)}(z,n)}{a(z)},Y_2^{(+)}(z,n)\right]\right),&z\in D_+,\\
		z^{n\hat\sigma_3}\left(\left[\begin{matrix}
			1&0\\0&c_n
		\end{matrix}\right]\left[Y_1^{(+)}(z,n),\frac{Y_2^{(-)}(z,n)}{\overline{a(\bar z^{-1}})}\right]\right),&z\in D_-.
	\end{cases}
\end{align}
Seeing from (\ref{e11s}), we have that
\begin{align}\label{e17ss}
	1-|r(z)|^2=\frac{c_{-\infty}}{|a(z)|^2},
\end{align}
and it follows that
\begin{align*}
	\parallel r\parallel_\infty<1. 
\end{align*} 
\begin{remark}\label{r2.2}
	$M$ possesses no pole on $\mathbb{C}\setminus\Sigma$ since $a(z)$ possesses no zero on $D_+$, i.e., there is no discrete spectrum. 
	Exactly, if there is a zero $z_0\in D_+$ for $a(z)$, then, by (\ref{e3}) and (\ref{e8s}), we see that there is a nonzero $2\times1$ $l^2$-vector function $f_0(n)$ that is an eigenfunction of (\ref{e2}):
	\begin{align}\label{e17s}
		f_0(n+1)=(z_0^{\sigma_3}+Q_n)f_0(n),
	\end{align}
	which is equivalent to
	\begin{align}\label{e18}
		c_{n}f_0(n)=(z_0^{\sigma_3}-Q_n)c_{n+1}f_0(n+1). 
	\end{align}
	Denoting $(f_0(n))^H$ as the Hermitian of $f_0(n)$, it's naturally to see the following equality
	\begin{align}
		\sum_{n=-\infty}^\infty (\sigma_3f_0(n+1))^H(c_nf_0(n))=\sum_{n=-\infty}^\infty (f_0(n))^H(\sigma_3c_{n-1}f_0(n-1)),
	\end{align}
	which deduces from (\ref{e17s}) and (\ref{e18}) that
	\begin{align}
		&\sum_{n=-\infty}^\infty (\tilde f_0(n))^HB_n\tilde f_0(n)=0\notag\\ 
		\Rightarrow& \sum_{n=-\infty}^\infty (\tilde f_0(n))^H\left(B_n+B_n^H\right)\tilde f_0(n)\notag\\
		&=(z_0+\bar z_0-z_0^{-1}-\bar z_0^{-1})\sum_{n=-\infty}^{\infty}(\tilde f_0(n))^H\tilde f_0(n)=0\notag\\
		\Rightarrow&z_0+\bar z_0-z_0^{-1}-\bar z_0^{-1}=0.\label{e20}
	\end{align}
	where
	\begin{align*}
		B_n=\sigma_3(z_0^{\sigma_3}-\bar z_0^{-\sigma_3}+Q_n-Q_{n-1}),\quad \tilde f_0(n)=\sqrt{c_n}f_0(n). 
	\end{align*}
	Similarly, we have
	\begin{align*}
		\sum_{n=-\infty}^\infty (f_0(n+1))^H(c_nf_0(n))=\sum_{n=-\infty}^\infty (f_0(n))^H(c_{n-1}f_0(n-1)),
	\end{align*}
	which, by the similar technique and the fact that $Q_n$ is Hermitian, implies that
	\begin{align}
		z_0-\bar z_0+z_0^{-1}-\bar z_0^{-1}=0.\label{e21}
	\end{align}
	Comparing (\ref{e20}) and (\ref{e21}), it's readily seen that $z_0\notin D_+$, i.e., $a(z)$ has no zero on $D_+$.
	As a result, $M$ is holomorphic on $\mathbb{C}\setminus\Sigma$.
\end{remark}
\begin{remark}\label{rm2.3}
	Recalling the fact that $a(z)$ is holomorphic and has no zero on $D_+$, by the relationship (\ref{e11s}) and (\ref{e9ss}), we obtain the trace formula for $a(z)$ and $z\in D_+$:
	\begin{align*}
		a(z)=\exp\left[-\frac{1}{2\pi i}\int_\Sigma\frac{\ln(1-|r(\zeta)|^2)}{\zeta-z}\ddddd\zeta\right].
	\end{align*}
\end{remark}
In the following, we discuss the RH problem about $M$. 
It's trivial to verify that $M$ admits RH problem \ref{r3}. The analyticity is obtained by Proposition \ref{l1},  (\ref{e9s}), Remark \ref{r2.2}. The normalization comes from (\ref{e10sb}) and (\ref{e9s}). And the jump condition comes from (\ref{e4}), (\ref{e8s}) and (\ref{e9s}). 
\begin{rhp}\label{r3}
	Look for a $2\times2$ function $M$ on $\mathbb{C}\setminus\Sigma$, such that:
	\begin{itemize}
		\item Analytisity:  $M$ is holomorphic on $\mathbb{C}\setminus\Sigma$. 
		\item Normalization:  as $z\to\infty$,
		\begin{align*}
			M\sim I+\mathcal{O}(z^{-1}). 
		\end{align*}
		\item Jump condition: on $z\in\Sigma$,
		\begin{align*}
			&M_+(z,n)=M_-(z,n)z^{n\hat\sigma_3}V(z),\\
			&V(z)=\left[\begin{matrix}
				1-|r(z)|^2&-\overline{r(z)}\\r(z)&1
			\end{matrix}\right].
		\end{align*}
		
	\end{itemize}
\end{rhp}

\begin{figure}
	\centering
	\includegraphics[width=0.3\linewidth]{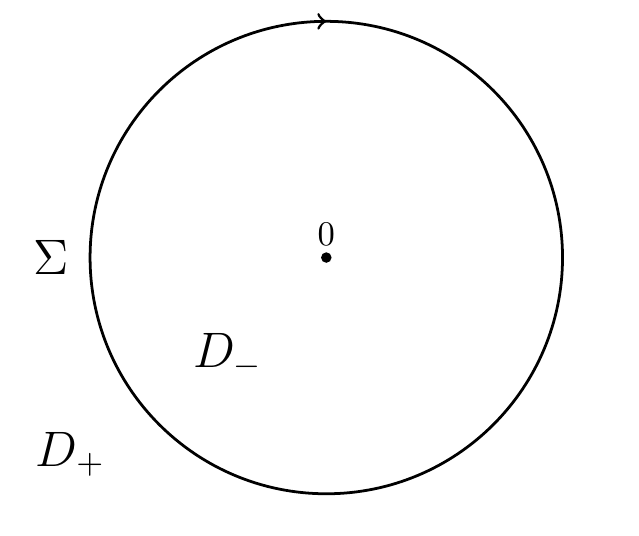}
	\caption{Jump contour: $\Sigma$, whose orientation is clockwise. $\mathbb{C}\setminus\Sigma=D_-\cup D_+$.}
\end{figure}

\subsection{Estimates of reflection coefficients}
\noindent

On the jump contour $\Sigma$, we write $z=e^{i\theta}$, $\theta\in[-\pi,\pi]$;
then, without confusion of notation, for $z\in\Sigma$, we write $r(\theta)=r(z)$ and  $Y^{(\pm)}(\theta,n)=Y^{(\pm)}(z,n)$. 
Then, we concentrate on some prior estimates and the proof of Theorem \ref{t2.4}.
\begin{theorem}\label{t2.4}
	Provided that potentials $q\in l^{2,k}$, reflection coefficients $r(\theta)\in H_\theta^k$. 
\end{theorem}
\begin{proof}
	By functional analysis, $H^k_\theta$ is a Banach algebra for $k\le1$. 
	If we denote $Y_{i,j}^{(\pm)}$ as the $i,j$-entry of $Y^{(\pm)}$, taking account of Proposition \ref{p2.6}, it follows that
	\begin{align}\label{e23}
		Y_{i,j}^{(\pm)}(\theta,0)\in H^k_\theta,\quad i,j=1,2. 
	\end{align} 
	From (\ref{e8}), (\ref{e23}) and the fact that $H^k_\theta$ is a Banach algebra, we see that
	\begin{align*}
		a(\theta), b(\theta)\in H_\theta^k.
	\end{align*}
	Considering (\ref{e11s}), we see that $a(\theta)\ne0$ for $\theta\in[-\pi,\pi]$, then it follows that $r(\theta)$ also belongs to $H^k_\theta$. 
\end{proof}
\noindent
In what follows, for proper notation, we simply denote all $L^2_\theta(\Sigma)$, $(L^2_\theta(\Sigma))^2$ and $(L^2_\theta(\Sigma))^{2\times2}$ as $L^2_\theta(\Sigma)$. 
The element in $L^2_\theta(\Sigma)$, $(L^2_\theta(\Sigma))^2$ and $(L^2_\theta(\Sigma))^{2\times2}$ are scalar, $2\times1$ vector and $2\times2$ matrix functions, respectively, and the $L^2$-norm are in general sense. 
We also denote
\begin{align*}
	\mathcal{A}=l^\infty(\mathbb{Z},L_\theta^2(\Sigma))=\left\{f=f(\theta,n):\sup_{n\in\mathbb{Z}}\parallel f(\cdot,n)\parallel_2<\infty\right\},
\end{align*}
and for any $f(\theta,n)\in\mathcal{A}$, the $\mathcal{A}$-norm denotes as
\begin{align*}
	\parallel f\parallel_{\mathcal{A}}=\sup_{n\in\mathbb{Z}}\parallel f(\cdot,n)\parallel_2.
\end{align*}
\begin{proposition}\label{p2.6}
	Provided that potentials $q\in l^{2,k}$, $k\ge1$, for any $\alpha=0,\dots,k$, we have $\partial_\theta^\alpha Y^{(\pm)}\in \mathcal{A}$.
\end{proposition}
\noindent For the proof of Proposition \ref{p2.6}, seeing (\ref{e5}), we write $Y_1^{(\pm)}(\theta,n)$ in the form of Neumann series:
\begin{align}\label{e9}
	Y_1^{(\pm)}(\theta,n)=\sum_{l=0}^{\infty}T_\pm^l\yl(\theta,n),\quad \theta\in[-\pi,\pi],
\end{align}
where $T_\pm$ are  operators
\begin{align*}
	T_- f(\theta,n)=e^{-i\theta\sigma_3}\sum_{k=-\infty}^{n-1}e^{-ik\theta\hat\sigma_3}Q_kf(\theta,k),\\
	T_+ f(\theta,n)=-e^{-i\theta\sigma_3}\sum_{k=n}^{\infty}e^{-ik\theta\hat\sigma_3}Q_kf(\theta,k).
\end{align*}
Then, it's readily seen that
\begin{subequations}\label{e10}
	\begin{align}
		&T_-^{2l-1}\yl(\theta,n)=\sum_{\scriptscriptstyle
			-\infty<k_{\scriptscriptstyle 2l-1}<\cdots< k_1< n
		}\left[\begin{matrix}
			0\\q_{k_2}\dots q_{k_{2l-2}}\overline{q_{k_1}\dots q_{k_{2l-1}}}e^{i\theta(A_{2l-1}+1)}
		\end{matrix}\right],\\
		&T_-^{2l}\yl(\theta,n)=\sum_{\scriptscriptstyle
			-\infty<k_{\scriptscriptstyle 2l}<\cdots< k_1< n
		}\left[\begin{matrix}
			q_{k_1}\dots q_{k_{2l-1}} \overline{q_{k_2}\dots q_{k_{2l}}}e^{-i\theta A_{2l}}\\0
		\end{matrix}\right],\label{e10b}\\
		&T_+^{2l-1}\yl(\theta,n)=\sum_{\scriptscriptstyle
			n\le k_1\le\cdots\le k_{\scriptscriptstyle 2l-1}< \infty
		}\left[\begin{matrix}
			0\\-q_{k_2}\dots q_{k_{2l-2}} \overline{q_{k_1}\dots q_{k_{2l-1}}}e^{i\theta(A_{2l-1}+1)}
		\end{matrix}\right],\\
		&T_+^{2l}\yl(\theta,n)=\sum_{\scriptscriptstyle
			n\le k_1\le\cdots\le k_{\scriptscriptstyle 2l}< \infty
		}\left[\begin{matrix}
			q_{k_1}\dots q_{k_{2l-1}} \overline{q_{k_2}\dots q_{k_{2l}}}e^{-i\theta A_{2l}}\\0
		\end{matrix}\right],
	\end{align}
\end{subequations}
where $l=1,2,\dots$, and 
\begin{align*}
	A_k=2\sum_{s=1}^{k}(-1)^{s-1}k_s. 
\end{align*}
\begin{lemma}\label{l2.7}
	Provided that the potential $q\in l^{2,k}$, $k\ge1$, we see that
	\begin{align*}
		\partial_\theta^\alpha T_\pm^l\yl\in \AAA,\quad l\ge1,\quad \alpha=0,\dots,k,
	\end{align*}
	and its $\mathcal{A}$-norm satisfies that
	\begin{subequations}\label{e24}
		\begin{align}
			&\Big\|\partial_\theta^\alpha T_\pm^{2l-1}\yl\Big\|_\AAA\le\left(1+\frac{2\pi^2}{3}\right)^{l-1}\frac{(4l-1)^\alpha\parallel q\parallel^{2l-1}_{2,k}}{(l-1)!^2},\\
			&\Big\|\partial_\theta^\alpha T_\pm^{2l}\yl\Big\|_\AAA\le\left(1+\frac{2\pi^2}{3}\right)^{l-\frac{1}{2}}\frac{(4l)^\alpha\parallel q\parallel^{2l}_{2,k}}{(l-1)!l!}.\label{e11b}
		\end{align}
	\end{subequations}
\end{lemma}
\begin{proof}
	To simplify the proof, we prove the case of $T_-^{2l}$, and other cases' proofs are similarly obtained. Ahead of all, we introduce an important result in functional analysis that for a scalar function $f(\theta)\in L_\theta^2(\Sigma)$,
	\begin{align}\label{e11}
		\parallel f\parallel_2=\sup_{g\in L_\theta^2(\Sigma),\parallel g\parallel_2<1}\Big|\frac{1}{2\pi i}\int_{-\pi}^{\pi}\overline{f(\theta)}g(\theta)\ddddd\theta\Big|. 
	\end{align}
	Defining a scalar function
	\begin{align*}
		&F(\theta,n)=\sum_{\scriptscriptstyle
			-\infty<k_{\scriptscriptstyle 2l}<\cdots< k_1< n
		}q_{k_1}\dots q_{k_{2l-1}} \overline{q_{k_2}\dots q_{k_{2l}}}e^{-i\theta(A_{2l})}
	\end{align*}
	that is the only nonzero entry of vector $T^{2l}_-\yl$, by (\ref{e11}) and the property of Fourier series, for $n\in\mathbb{Z}$, we obtain estimates of the $L^2$-norm 
	\begin{align}
		\parallel F(\cdot,n)\parallel_2&=\sup_{g\in L_\theta^2(\Sigma),\parallel g\parallel_2<1}\sum_{\scriptscriptstyle
			-\infty<k_{\scriptscriptstyle 2l-1}<\cdots< k_1< n
		}|q_{k_1}\dots q_{k_{2l-1}} \overline{q_{k_2}\dots q_{k_{2l}}}\hat g(A_{2l})|\notag\\
		&\le\parallel q\parallel_2\sum_{\scriptscriptstyle
			-\infty<k_{\scriptscriptstyle 2l-1}<\cdots< k_1< n
		}|q_{k_1}\dots q_{k_{2l-1}} \overline{q_{k_2}\dots q_{k_{2l-2}}}|\notag\\
		&\le\frac{\parallel q\parallel_1^{2l-1}\parallel q\parallel_2}{l!(l-1)!}\le \left(1+\frac{2\pi^2}{3}\right)^{l-\frac{1}{2}}\frac{\parallel q\parallel_{2,k}^{2l}}{l!(l-1)!},\label{e12}
	\end{align}
	where the last inequality holds by Schwartz inequality
	\begin{align*}
		\parallel q\parallel_1\le\sqrt{\sum_{n=-\infty}^{\infty}\frac{1}{1+n^2}}\parallel q\parallel_{2,1}\le\sqrt{1+2\sum_{1}^{\infty}\frac{1}{n^2}}\parallel q\parallel_{2,1}\le\sqrt{1+\frac{2\pi^2}{3}}\parallel q\parallel_{2,k}.
	\end{align*}
	As a result of (\ref{e12}), we obtain that
	\begin{align*}
		\parallel F\parallel_\mathcal{A}\le\left(1+\frac{2\pi^2}{3}\right)^{l-\frac{1}{2}}\frac{\parallel q\parallel_{2,k}^{2l}}{l!(l-1)!}.
	\end{align*}
	Differentiate $F(\theta,n)$ about $\theta$:
	\begin{align}\label{e15}
		\partial_\theta F(\theta,n)=-2i\sum_{\tilde l=1}^{2l}(-1)^{\tilde l-1}F_{\tilde l}(\theta,n),
	\end{align}
	where for $\tilde l=1,\dots,2l$,
	\begin{align*}
		F_{\tilde l}(\theta,n)=\sum_{\scriptscriptstyle
			-\infty<k_{\scriptscriptstyle 2l-1}<\cdots< k_1< n
		}k_{\tilde l}q_{k_1}\dots q_{k_{2l-1}} \overline{q_{k_2}\dots q_{k_{2l}}}e^{-i\theta(A_{2l})}. 
	\end{align*}
	Again using (\ref{e11}), properties of Fourier series and the Schwartz inequality, it is readily seen that for all $\tilde l=1,\dots,2l$, 
	\begin{align*}
		\parallel F_{\tilde l}\parallel_\mathcal{A}\le \left(1+\frac{2\pi^2}{3}\right)^{l-\frac{1}{2}}\frac{\parallel q\parallel_{2,k}^{2l}}{l!(l-1)!},
	\end{align*}
	which implies by (\ref{e15}) that
	\begin{align}\label{e16}
		\parallel \partial_\theta F\parallel_\mathcal{A}\le4l\left(1+\frac{2\pi^2}{3}\right)^{l-\frac{1}{2}}\frac{\parallel q\parallel_{2,k}^{2l}}{l!(l-1)!}. 
	\end{align}
	Repeat these procedure, and then it naturally follows that for $\alpha=0,\dots,k$,
	\begin{align}\label{e17}
		\parallel\partial_\theta^\alpha F\parallel_\mathcal{A}\le(4l)^\alpha\left(1+\frac{2\pi^2}{3}\right)^{l-\frac{1}{2}}\frac{\parallel q\parallel_{2,k}^{2l}}{l!(l-1)!}.
	\end{align}
	Since $F(\theta,n)$ is the only non-vanishing entry of the vector $T_-^{2l}\yl(\theta,n)$, we see from (\ref{e17}) that $T_-^{2l}\yl$ satisfies (\ref{e11b}) for $\alpha=0,\dots,k$. 
	The remaining result for $T_-^{2l-1}$, $T_+^{2l-1}$ and $T_+^{2l}$ in (\ref{e24}) can be similarly obtained. The result is valid.
\end{proof}

\begin{proof}[Proof of Proposition \ref{p2.6}]
	Here, we only prove the part of $\partial_\theta^\alpha Y_1^{(\pm)}$, and the proof about $\partial_\theta^\alpha Y_2^{(\pm)}$ is parallel. (\ref{e24}) implies that
	\begin{align}\label{e30}
		&\Big\|\partial_\theta^\alpha T_\pm^l\yl\Big\|_\AAA\le a_{l,\alpha}\parallel q\parallel_{2,k}^l, \quad
		&a_{l,\alpha}=\begin{cases}
			\frac{(1+\frac{2\pi^2}{3})^{\frac{l-1}{2}}}{(\frac{l-1}{2})!^2}(2l+1)^\alpha,&\text{l is odd},\\
			\frac{(1+\frac{2\pi^2}{3})^{\frac{l-1}{2}}}{(\frac{l-2}{2})!(\frac{l}{2})!}(2l)^\alpha,&\text{l is even}.
		\end{cases}
	\end{align}
	By trivial computation, we see from (\ref{e30}) that
	\begin{align*}
		\lim_{l\to\infty}\frac{a_{l+1,\alpha}}{a_{l,\alpha}}=0,
	\end{align*}
	and then it follows that the infinite summation
	\begin{align*}
		\sum_{l=0}^{\infty}a_{l,\alpha}\parallel q\parallel_{2,k}^l
	\end{align*}
	converges for all $q\in l^{2,k}$, $\alpha=0,\dots,k$; therefore, the result in Proposition \ref{p2.6} is naturally confirmed by control convergence theorem, Lemma \ref{l2.7}, (\ref{e9}) and (\ref{e30}). 
\end{proof}

\section{Inverse scattering problem}\label{s3}
\indent

In Section \ref{s2}, we have proven that if potentials satisfy (\ref{e3s}), then we have $r\in H^k_\theta$ and $\parallel r\parallel_\infty<1$. 
In this section,  we consider inverse scattering mapping from the reflection coefficient $r$ to the potential $q$ by obtaining the reconstructed formula for $q$, and the purpose is to check Theorem \ref{th1}. 
\begin{theorem}\label{th1}
	Given $r\in H^k_\theta$ and $\parallel r\parallel_\infty<1$, the potential $q$ belongs to $l^{2,k}$.
\end{theorem}

\subsection{Reconstructed formulas}
\noindent

Here, we obtain reconstructed formulas, and base on them, we recover the potential for $n\in\mathbb{Z}$. 
Introducing $2\times2$ matrix functions:
\begin{align*}
	&w=w_++w_-,\quad 
	w_-\equiv w_-(\theta,n)\equiv w_-(z,n)=\left[\begin{matrix}
		0&-\overline{r(z)}z^{2n}\\0&0
	\end{matrix}\right],\\
&w_+\equiv w_+(\theta,n)\equiv w_+(z,n)=\left[\begin{matrix}
	0&0\\r(z)z^{-2n}&0
\end{matrix}\right],
\end{align*}
where $z=e^{i\theta}\in \Sigma$, 
it's readily seen that
\begin{align*}
	z^{n\hat\sigma_3}V(z)=(I-w_-(z,n))^{-1}(I+w_+(z,n)).
\end{align*}
Also, introducing a Cauchy-type operator
\begin{align}
	&C_wf=C_+(fw_-)+C_-(fw_+),\label{e32s}\\
	&C_\pm f(z)=\lim_{z'=(1+\epsilon)z, \epsilon\to\pm0}\frac{1}{2\pi i}\int_\Sigma\frac{f(z')}{z'-z}\mathrm{d}z',
\end{align}
it's well known from \cite{ablowitz2003complex} that the Cauchy integral operator $C_\pm$ act on $L^2_\theta(\Sigma)$ and
\begin{align}\label{e29}
	\parallel C_\pm\parallel_{L^2\to L^2}\le1,\quad C_+-C_-=I. 
\end{align}
Considering the Beals-Coifman solution for RH problem \ref{r3}
\begin{align}\label{e35ss}
	M(z,n)=I+\frac{1}{2\pi i}\int_\Sigma\frac{[(I-C_w)^{-1}I](z',n)w(z',n)\ddddd z'}{z'-z}, 
\end{align}
we derive the reconstructed formula by (\ref{e10sa}), (\ref{e9s}) and (\ref{e35ss}):
\begin{align}\label{e31s}
q_n&=\lim_{z\to0}\frac{[M(z,n+1)]_{1,2}}{z}=\partial_z[M(z,n+1)]_{1,2}\big|_{z=0}\notag\\
&=\frac{1}{2\pi i}\int_\Sigma z^{-2}\left[(I-C_w)^{-1}Iw\right]_{1,2}(z,n+1)\mathrm{d}z. 
\end{align}
Define a new matrix function
\begin{align}\label{e37s}
	\tilde M=(c_{-\infty})^{-\frac{\hat\sigma_3}{2}}M\delta^{-\sigma_3},
\end{align}
where $\delta\equiv\delta(z)=\exp\left[\frac{1}{2\pi i}\int_\Sigma\frac{\ln(1-|r(\zeta)|^2)}{\zeta-z}\mathrm{d}\zeta\right]$ is the unique solution of RH problem \ref{r3.2} and $c_{-\infty}$ is defined as (\ref{e6}), and we obtain another reconstructed formula (\ref{e42s}) about this RH problem in the following.
\begin{remark}
	Seeing from the trace formula in Remark \ref{rm2.3}, we can recover $a(z)$ from the reflection coefficient $r(z)$. 
	Moreover, $a(z)$ is naturally holomorphic and has no zero on $D_+$. 
	Then, we discuss the equivalence of two definitions of $c_{-\infty}$ shown in Introduction.
	Considering condition (\ref{e17ss}), we find that
	\begin{align}\label{e36s}
		\frac{1}{2\pi i}\frac{\ln(1-|r(z)|^2)}{z}=\frac{1}{2\pi i}\frac{\ln c_{-\infty}-\ln |a(z)|^2}{z},\quad z\in\Sigma.
	\end{align}
	Integrating both sides of (\ref{e36s}) and then exponentiating it, we learn that
	\begin{align}\label{e37}
		c_{-\infty}=\delta(0)\exp\left[\frac{1}{2\pi}\int_{-\pi}^{\pi}\ln|a(\theta)|^2\ddddd\theta\right].
	\end{align}
	Recalling the fact that $a(z)$ is holomorphic and has no zero on $D_+$, we have that
	\begin{align}\label{e35s}
		\frac{1}{2\pi}\int_{-\pi}^{\pi}\ln|a(\theta)|^2\ddddd\theta&=\frac{1}{2\pi i}\int_\Sigma\frac{\ln a(z)}{z}\ddddd z+\frac{1}{2\pi i}\int_\Sigma\frac{\ln \overline{a(\bar z^{-1})}}{z}\ddddd z\notag\\
		&=\lim_{\epsilon\to+0}\left(\frac{1}{2\pi i}\int_{\epsilon^{-1}\Sigma}\frac{\ln a(z)}{z}\ddddd z+\frac{1}{2\pi i}\int_{\epsilon\Sigma}\frac{\ln \overline{a(\bar z^{-1})}}{z}\ddddd z\right)\notag\\
		&=\lim_{z\to\infty}(\ln a(z)+\ln \overline{a(z)})=0. 
	\end{align}
Comparing (\ref{e37}) and (\ref{e35s}), we confirm that the definition of $c_{-\infty}$ in (\ref{e6}) is equivalent to that in (\ref{e5s}). 
\end{remark}
\begin{rhp}\label{r3.2}
	Find a scalar function $\delta$ such that
	\begin{itemize}
		\item Analyticity: $\delta$ is holomorphic on $\mathbb{C}\setminus\Sigma$;
		\item Normalization: as $z\to\infty$, $\delta(z)\to 1$;
		\item Jump condition: on $\Sigma$, we have $\delta_+=\delta_-(1-|r|^2)$. 
	\end{itemize}
\end{rhp}
\noindent
By (\ref{e6}), (\ref{e37s}) and RH problem \ref{r3.2}, it's readily seen that RH problem \ref{r3} is equivalent to  RH problem \ref{r4}:
\begin{rhp}\label{r4}
	Find a $2\times2$ matrix function $\tilde M\equiv \tilde M(z,n)$ such that,
	\begin{itemize}
		\item Analyticity: $\tilde M$ is holomorphic on $\mathbb{C}\setminus\Sigma$;
		\item Normalization: as $z\to\infty$,
		\begin{align*}
			\tilde M(z,n)\to I;
		\end{align*}
		\item Jump condition: on $z\in\Sigma$,
		\begin{align*}
			&\tilde M_+(z,n)=\tilde M_-(z,n)z^{n\hat\sigma_3}\tilde V(z),\\
			&\tilde V(z)=\left[\begin{matrix}
				1&-\overline{\tilde r(z)}\\\tilde r(z)&1-|\tilde r(z)|^2
			\end{matrix}\right],
		\end{align*}
		where
		\begin{align*}
			\tilde r(z)=r(z)\delta_-^{-1}(z)\delta_+^{-1}(z)c_{-\infty},\quad z\in\Sigma. 
		\end{align*}
	\end{itemize}
	By (\ref{e37s}) and RH problem \ref{r4}, we similarly obtain the reconstructed formula with respect to $\tilde w$
	\begin{align}\label{e42s}
		q_n=\frac{1}{2\pi i}\int_{\Sigma}z^{-2}[(I-C_{\tilde w})^{-1}I\tilde w]_{1,2}(z,n+1)\ddddd z, 
	\end{align}
	which is equivalent to (\ref{e31s}), where 
	\begin{align*}
		&\tilde w=\tilde w_-+\tilde w_+,\quad
		\tilde w_-\equiv\tilde w_-(z,n)=\left[\begin{matrix}
			0&0\\\tilde r(z)z^{-2n}&0
		\end{matrix}\right],\\
	 &\tilde w_+\equiv\tilde w_+(z,n)=\left[\begin{matrix}
			0&-\overline{\tilde r(z)}z^{2n}\\0&0
		\end{matrix}\right].
	\end{align*}
	By basic computation, we have that $z^{n\hat\sigma_3}\tilde V(z)$ admits a factorization 
	\begin{align*}
		z^{n\hat\sigma_3}\tilde V(z)=(I-\tilde w_-(z,n))^{-1}(I+\tilde w_+(z,n)). 
	\end{align*}
\end{rhp}

\subsection{Estimates of potentials}
\noindent

Here, we aim to completing proof of Theorem \ref{th1}.
By definition (\ref{e32s}), we also introduce the operator on $L^2_\theta(\Sigma)$, $C_{\tilde w}$. 
For operators $C_w$ and $C_{\tilde w}$, we make some prior estimates associated to them and obtain the following lemmas. 
\begin{lemma}\label{l2}
	When $\parallel r\parallel_\infty<1$, $1-C_w$ and $1-C_{\tilde w}$ are both invertible on $L_\theta^2(\Sigma)$. Moreover, $(1-C_w)^{-1}$ and $(1-C_{\tilde w})^{-1}$ are both bounded on $L_\theta^2(\Sigma)$. 
\end{lemma}
\begin{lemma}\label{l3}
	When $r(\theta)\in H^k_\theta$ and $n\ge0$, $C_wI$ belongs to $L_\theta^2(\Sigma)$ and for some constant $C>0$,
	\begin{align*}
		\parallel C_wI\parallel_2\le C (1+n^2)^{-\frac{k}{2}}.
	\end{align*}
\end{lemma}
\begin{lemma}\label{l4}
	When $r(\theta)\in H^k_\theta$, $\parallel r\parallel_\infty<1$ and $n<0$, $C_{\tilde w}I$ belongs to $L_\theta^2(\Sigma)$ and  for some constant $C>0$,
	\begin{align*}
		\parallel C_{\tilde w}I\parallel_2\le C (1+n^2)^{-\frac{k}{2}}.
	\end{align*}
\end{lemma}
\begin{proof}[proof of Lemma \ref{l2}]
	We first check the part of $C_w$. 
	For any $2\times2$ matrix function $f(\theta)\in L_\theta^2(\Sigma)$, we claim that
	\begin{align}\label{e31}
		C_wf=\left[C_+(\bar re^{2in\cdot}f_2),C_-(re^{-2in\cdot}f_1)\right],\quad f=[f_1,f_2]. 
	\end{align}
It follows by (\ref{e29}) and (\ref{e31}) that
\begin{align*}
	\parallel C_wf\parallel_2\le \parallel r\parallel_\infty\parallel f\parallel_2. 
\end{align*}
and therefore, we claim that $1-C_w$ is invertible and $(1-C_w)^{-1}$ is bounded on $L_\theta^2(\Sigma)$ with the fact that $\parallel r\parallel_\infty<1$. 

Then, we come to the part of $C_{\tilde w}$. 
Since $\delta(z)$ solves RH problem \ref{r3.2}, by (\ref{e6}), we learn that $c_{-\infty}\left(\overline{\delta(\bar z^{-1})}\right)^{-1}$ is also the solution of RH problem \ref{r3.2} and by the uniqueness of solution, we have
\begin{align*}
	\delta(z)=c_{-\infty}\left(\overline{\delta(\bar z^{-1})}\right)^{-1},\quad z\in\mathbb{C}\setminus\Sigma,
\end{align*}
which implies that
\begin{align}\label{e41}
	c_{-\infty}\delta_+^{-1}(z)\delta_-^{-1}(z)=\overline{\delta_-(z)}\delta_-^{-1}(z),\quad z\in\Sigma.
\end{align}
It follows from (\ref{e41}) that
\begin{align}\label{e42}
	\parallel \tilde r\parallel_\infty=\parallel r\parallel_\infty<1.
\end{align}
Therefore, the proof for $(1-C_{\tilde w})^{-1}$ is parallel to that of $(1-C_{w})^{-1}$ and the result is valid.
\end{proof}
\begin{proof}[proof of Lemma \ref{l3}]
	By simple computation, it's readily seen that
	\begin{align}\label{e32}
		C_wI=C_+w_-+C_-w_+,\quad C_-w_+=\left[\begin{matrix}
			0&r_1\\0&0
		\end{matrix}\right], \quad C_+w_-=\left[\begin{matrix}
		0&0\\r_2&0
	\end{matrix}\right], 
	\end{align}
where
\begin{align}\label{e33}
	r_1(\theta,n)=\sum_{l=2n+1}^{\infty}-\overline{\hat r(l)}e^{i(2n-l)\theta},\quad r_2(\theta,n)=\sum_{l=2n}^\infty\hat r(l)e^{i(l-2n)\theta}.
\end{align}
By (\ref{e33}), Parseval theorem and Schwartz Inequality,  it follows that
\begin{align*}
	&\parallel r_1(\cdot,n)\parallel_2=\left(\sum_{l=2n+1}^{\infty}|\hat r(l)|^2\right)^\frac{1}{2}\le (1+n^2)^{-\frac{k}{2}}\parallel r\parallel_{H^k_\theta},\\
	&\parallel r_2(\cdot,n)\parallel_2=\left(\sum_{l=2n}^{\infty}|\hat r(l)|^2\right)^\frac{1}{2}\le (1+n^2)^{-\frac{k}{2}}\parallel r\parallel_{H^k_\theta},
\end{align*}
which combined with (\ref{e32}) confirms the result. 
\end{proof}
\begin{proof}[proof of Lemma \ref{l4}]
	Since $r\in H^k_\theta$ and $\parallel r\parallel_\infty<1$, it's readily seen that
	\begin{align*}
		\rho\deff\ln(1-|r|^2)\in H^k_\theta.
	\end{align*}
	By Fourier theory, we learn that
	\begin{align*}
		C_+\rho=\sum_{l=-\infty}^{-1}\hat\rho(l)e^{il\theta},\quad C_-\rho=-\sum_{l=0}^{\infty}\hat\rho(l)e^{il\theta}.
	\end{align*}
	and they  both belong to $H^k_\theta$. 
	Since it's trivial to check that
	\begin{align*}
		\delta_+=\exp[C_+\rho], \quad \delta_-=\exp[C_-\rho],
	\end{align*}
	and $H^k_\theta$ is a Banach algebra, we also have that 
	\begin{align*}
		\delta_+,\delta_-\in H^k_\theta.
	\end{align*}
	Since $r,\delta_+,\delta_-\in H^k_\theta$ and $H^k_\theta$ is a Banach algebra, by the definition of $\tilde r$, we claim that
	\begin{align}\label{e40}
		\tilde r\in H^k_\theta.
	\end{align}
	Because of (\ref{e42}) and (\ref{e40}), the remained proof is parallel to that of Lemma \ref{l3}.
\end{proof}
\begin{proof}[proof  of Theorem \ref{th1}]
	To prove Theorem \ref{th1}, We firstly prove that
	\begin{align}\label{e49}
		\sum_{n=-1}^{\infty}(1+n^2)^k|q_n|^2<\infty,
	\end{align}
	and secondly prove that
	\begin{align}\label{e50}
		\sum_{n=-\infty}^{-2}(1+n^2)^k|q_n|^2<\infty.
	\end{align}

	For (\ref{e49}), from (\ref{e31s}), we derive that
	\begin{align}\label{e35}
		q_{n-1}&=\frac{1}{2\pi}\int_{-\pi}^{\pi}-\overline {r(\theta)}e^{2in\theta}+\left[C_wIw\right]_{1,2}(\theta,n)+[(I-C_w)^{-1}(C_w^2I)w]_{1,2}(\theta,n)\mathrm{d}\theta\notag\\&=I_1+I_2+I_3.
	\end{align}
Since $r\in H^k_\theta$, $I_1$ belongs to $l^{2,k}$ by Fourier analysis, and $I_2=0$ by trivial computation. 
Then, we discuss the part about $I_3$. 
Setting
\begin{align*}
	\mu=(I-C_w)^{-1}(C_wI)
\end{align*}
by Cauchy's theorem and (\ref{e29}), we have
\begin{align}\label{e36}
	I_3&=\frac{1}{2\pi}\int_{-\pi}^{\pi}[(C_w\mu)w]_{1,2}(\theta,n)\ddddd\theta\notag\\
	&=\frac{1}{2\pi}\int_{-\pi}^{\pi}[C_+(\mu w_-)C_-w_+]_{1,2}(\theta,n)+[C_-(\mu w_+)C_+w_-]_{1,2}(\theta,n)\mathrm{d}\theta
\end{align}
By Schwartz inequality, (\ref{e29}), (\ref{e36}) Lemma \ref{l2} and Lemma \ref{l3},
\begin{align*}
	|I_3|&\le \parallel(\mu w_-)(\cdot,n)\parallel_2\parallel C_-w_+(\cdot,n)\parallel_2+\parallel(mu w_+)(\cdot,n)\parallel_2\parallel C_+w_-(\cdot,n)\parallel_2\\
	&\le\parallel\mu(\cdot,n)\parallel_2\parallel w(\cdot,n)\parallel_\infty\parallel C_wI(\cdot,n)\parallel_2\\
	&\le\parallel r\parallel_\infty\parallel (I-C_w)^{-1}\parallel_{L^2\to L^2}\parallel C_wI(\cdot,n)\parallel_2^2\le C(1+n^2)^{-k},
\end{align*}
where $C$ is only dependent on $\parallel r\parallel_{H^k_\theta}$. Therefore, we obtain (\ref{e49}). 

For (\ref{e50}), because of Lemma \ref{l2} and \ref{l4}, this part is parallel to the part of (\ref{e49}). The result is confirmed when taking account of (\ref{e49}) and (\ref{e50}). 
\end{proof}

\section{Time evolution}\label{s4}
\indent

In this section, we consider the solution of defocusing Ablowitz-Ladik systems (\ref{e1}) with initial potential $q(0)=\{q_n(0)\}_{-\infty}^\infty$ satisfying (\ref{e3s}). 
The Lax-pair for defocusing Ablowitz-Ladik systems \cite{Ablowitz1975nonlinear} is 
\begin{align*}
	V(z,n+1,t)=&\left[\begin{matrix}
		z&q_n(t)\\\overline{q_n(t)}&1/z
	\end{matrix}\right]V(z,n,t),\\
	\partial_tV(z,n,t)=&\left[\begin{matrix}
		-i(\frac{1}{2}(z-1/z)^2-q_n(t)\overline{q_{n-1}(t)})&(q_{n-1}(t)/z-q_n(t)z)\\
		i(\overline{q_n(t)}/z-\overline{q_{n-1}(t)}z)&i(\frac{1}{2}(z-1/z)^2-\overline{q_n(t)}q_{n-1}(t))
	\end{matrix}\right]\\&\times V(z,n,t),
\end{align*}
where the spatial part is exactly the spectral problem (\ref{e2}). 

We denote $r(\theta,t)=r(z,t)$ as the reflection coefficient for (\ref{e1}) when $t\in \mathbb{R}$, where $r(\cdot,0)$ is the reflection coefficient for the initial potential and 
\begin{align*}
	z=e^{i\theta}\in\Sigma
\end{align*}
From section \ref{s2}, it's readily seen that if $q(0)$ satisfying (\ref{e3s}), then $r(\cdot,0)\in H^k_\theta$ and $\parallel r(\cdot,0)\parallel_\infty<1$. 
Then, we study how the reflection coefficient evolves about time parameter $t$. 
Seeing Section 3 in \cite{ablowitz2004discrete}, it is readily seen that
\begin{align}\label{e38}
	r(\theta,t)=r(\theta,0)e^{2i(\cos2\theta-1)t},\quad \theta\in[-\pi,\pi].
\end{align}
and $\parallel r(\cdot,t)\parallel_\infty<1$ for all fixed $t\ne0$.
It's easy to check that for any $k\ge1$,
\begin{align*}
	\phi(\theta)=\cos2\theta-1\in H^k_\theta, \quad \parallel\phi\parallel_{H^k_\theta}=1+\frac{5^k}{2}.
\end{align*}
Therefore, since $H^k_\theta$ is a Banach algebra, it follows that
\begin{align*}
	e^{2it\phi(\theta)}\in H^k_\theta, \quad \text{and}\quad \parallel e^{2it\phi}\parallel_{H^k_\theta}\le e^{2t(1+\frac{5^k}{2})}. 
\end{align*}
which combined with (\ref{e38}) deduces that
\begin{align*}
	r(\cdot,t)\in H^k_\theta.
\end{align*}
It follows from Section \ref{s3} that we can recover $q_n(t)$ by the reconstructed formula, and that $q_n(t)$ still belongs to $l^{2,k}$, i.e., we solve (\ref{e1}) for the initial data satisfying (\ref{e3s}).

\section*{Acknowledgements}
\noindent

This work is supported by the National Natural Science of China (Grant No. 12071304). 
The authors would like to have their sincerest gratitude to referees for patient guidance and valuable suggestions. 

\section*{Declaration of competing interests}
\noindent

The authors declare that they have no conflict of interest.

\end{document}